# ROBUSTNESS OF MULTIPLE TESTING PROCEDURES AGAINST DEPENDENCE

By Sandy Clarke and Peter Hall

*University of Melbourne*

An important aspect of multiple hypothesis testing is controlling the significance level, or the level of Type I error. When the test statistics are not independent it can be particularly challenging to deal with this problem, without resorting to very conservative procedures. In this paper we show that, in the context of contemporary multiple testing problems, where the number of tests is often very large, the difficulties caused by dependence are less serious than in classical cases. This is particularly true when the null distributions of test statistics are relatively light-tailed, for example, when they can be based on Normal or Student's $t$ approximations. There, if the test statistics can fairly be viewed as being generated by a linear process, an analysis founded on the incorrect assumption of independence is asymptotically correct as the number of hypotheses diverges. In particular, the point process representing the null distribution of the indices at which statistically significant test results occur is approximately Poisson, just as in the case of independence. The Poisson process also has the same mean as in the independence case, and of course exhibits no clustering of false discoveries. However, this result can fail if the null distributions are particularly heavy-tailed. There clusters of statistically significant results can occur, even when the null hypothesis is correct. We give an intuitive explanation for these disparate properties in light- and heavy-tailed cases, and provide rigorous theory underpinning the intuition.

**1. Introduction.** Classical properties of simultaneous hypothesis testing, error rate and false-discovery rate are well understood. They have been explored extensively, in both practice and theory, in the context of independent tests. However, for a range of contemporary applications, multiple testing

---

problems differ substantially from the conventional. For instance, the number, $\nu$ say, of tests is often far greater than the number of data, $n$, in the samples from which test statistics are computed. There is also potential for a degree of dependence among samples, even though the data within a sample can often fairly be assumed to be independent.

By way of contrast, in classical settings the value of $\nu$ is relatively small, and critical points are only moderately large (equivalently, $p$-values are only modestly small). Here a major, noticeable impact of dependence is that it results in clusters of rejections. That is, if a test is rejected for a particular value of an index, then there are likely to be further rejections for tests that have nearby indices (assuming that index order reflects dependence). This can impact significantly on the accuracy of multiple testing procedures.

One approach to alleviating the difficulties caused by dependence is to use techniques based on Bonferroni's inequality. However, such bounds are quite conservative, and if they could be avoided, then greater precision would result. In some settings, where positive dependence is present, corrections of Bonferroni type are unnecessary [see Benjamini and Yekutieli (2001) for discussion], but in general the nature of dependence is not known reliably. Moreover, even in the case of positive dependence it is of interest to know whether the test is genuinely conservative, as indicated by conventional theoretical arguments, or whether its level accuracy is virtually the same as in the case of independent data. Efron (2007) has suggested correlation corrections for large-scale simultaneous hypothesis testing.

One might expect the same difficulties and questions to arise in contemporary testing problems, where $\nu$ is much greater than $n$. (In some of these problems, typical values of $\nu$ and $n$ are 10,000 and 20, resp.) Indeed, there is reason to suspect that difficulties could increase with increasing $\nu$, since it can be particularly difficult to model accurately the extremes of dependent data processes. Additionally, inaccuracies become more obvious as the amount of information about a model increases.

However, it turns out that sometimes, although not always, the problem is actually simpler in the contemporary, "$\nu$ much larger than $n$" case. For example, in cases where test statistics have light-tailed distributions, the difficulties caused by dependence tend to retreat as the number of simultaneous tests increases. The number of clusters of false discoveries declines, and the distribution of critical-point exceedances closely resembles its counterpart for independent data. Only for very heavy-tailed data is this property violated; for dependent data, when the distribution of the test statistic is light-tailed and the number of simultaneous tests is very large, methods that would normally be recommended only for independent data can give good control of error rate and false-discovery rate.

This result can be explained intuitively by noting that, in the case of light-tailed marginal distributions, exceedances above a high level occur only because neighboring disturbances are fortuitously aligned. Indeed, since the tail is light, then it is highly unlikely that a single disturbance is so great as to carry the process close to, or over, the level for several different indices. Instead different, moderately large disturbances reinforce one another, by chance, at a particular index. However, at adjacent indices the circumstances that led to alignment change. As a result the propensity for level exceedence quickly diminishes, and even disappears. Consequently, clusters of exceedences seldom arise. That is, the pattern of exceedences appears as though it was produced by a sequence of independent tests, and as a result, both generalized family-wise error rate, and false-discovery rate, can be controlled by appealling to standard arguments for independent tests.

On the other hand, when test statistics have heavy-tailed distributions it is possible for a single disturbance to be so great that it carries the value of a test statistic over a high level for several indices in a row. In such cases, clusters of exceedences occur, and methods based on independent data are not adequate for controlling error rates.

These arguments and properties, especially those in the light-tailed setting, are applicable only to exceedences of high levels. Very high levels are relevant only when the number of simultaneous tests is very large, and so the properties tend not to be noticed in conventional multiple testing problems, where the number of tests is relatively small.

In this paper we develop rigorous arguments, using linear-process models for test statistics, to capture in theory the ideas discussed above. We show that if the test statistic distribution has tails that decay like $\exp(-Cx^\gamma)$, for constants $C, \gamma > 0$, then the tails can be regarded as "light" (in the context of the discussion above) when $\gamma \geq 1$; they are "heavy" when $0 < \gamma < 1$. However, even in the latter case the problem has many of the characteristics of the light-tailed context, unless there are ties among the weights in the linear process. Only in very heavy-tailed cases, where the distribution of the test statistic decreases at a polynomial rather than exponential rate in the tails, are methods based on independent data seen to be inadequate.

Moreover, even in these heavy-tailed contexts the independent-data approach can provide good results for large-but-not-too-large $\nu$. A case in point is that where the test statistic is a Student's $t$ ratio. There, although the extreme tails of the test statistic distribution are typically regularly varying (e.g., when the sampling distribution is Gaussian), large-deviation properties show that less extreme parts of the tail are well approximated by the function $\exp(-Cx^\gamma)$ for $\gamma = 2$ [see, e.g., Shao (1999) and Wang (2005)]. As a result, good performance can be obtained, in the case of dependent $t$-statistics, by arguing as though the data are independent.

There is a particularly broad and deep literature on multiple testing procedures, only a part of it confined to statistics journals. Review-type contributions include those of Hochberg and Tamhane (1987), who expounded work on multiple comparisons up to the mid-1980s; Pigeot (2000), who surveyed conceptual issues in multiple testing; Dudoit, Shaffer and Boldrick (2003), who reviewed multiple hypothesis testing in microarray settings; Bernhard, Klein and Hommel (2004), who discussed literature on global and multiple testing; and [Lehmann and Romano (2005), Chapter 9], who discussed multiple hypothesis testing in the context of hypothesis testing more generally.

Among contributions related to this paper, Hochberg and Benjamini (1990) pointed to the need for procedures that are more powerful than classical multiple comparison methods, and suggested new, generally applicable techniques; Rom (1990) introduced methods based on modified Bonferroni arguments; Dunnett and Tamhane (1995) discussed step-up methods for multiple testing in the presence of correlation; Wright (1992) developed $p$-value adjustments based on Bonferroni's bounds; Benjamini and Hochberg (2000, 1995) proposed approaches to false-discovery rate in multiple testing; Blair, Troendle and Beck (1996) introduced methods for controlling family-wise error rates in multiple procedures; Brown and Russell (1997) suggested corrections for multiple testing; Olejnik et al. (1997) compared Bonferroni-type methods; Sarkar and Chang (1997) discussed multiple testing in the presence of positive dependence; Finner and Roters (1998, 1999, 2000) gave asymptotic theory for an increasingly large number of hypothesis tests, and (2002) discussed the expected number of Type I errors in multiple testing problems; Holland and Cheung (2002) discussed robustness of family-wise error rate; Kesselman, Cribbie and Holland (2002) suggested ways of controlling level accuracy over a large number of hypothesis tests; Genovese and Wasserman (2004) proposed new, stochastic process-based methods for controlling false-discovery rate in multiple testing; Lehmann, Romano and Shaffer (2005) developed optimality theory for multiple testing; Rosenberg, Che and Chen (2006) suggested multiple hypothesis testing methods in a genomic setting; Sarkar (2006) obtained new results on false-discovery rates for single-step, multiple testing procedures; Schmidt and Stadtmüller (2006) and Schmidt (2007) discussed upper-tailed dependence; and Yekutieli et al. (2006) developed new approaches to the treatment of multiplicity in the setting of microarray analysis.

The issue of overall error rate, as distinct from the error rate of individual tests, was taken up by Godfrey (1985), who drew attention to the tendency to enhance the significance of treatment effects if the overall error rate is not controlled. See also [Smith et al. (1987), Pocock, Hughes and Lee (1987), Gotzsche (1989), Ludbrook (1991), Ottenbacher (1991a, 1991b, 1998) and Ottenbacher and Barrett (1991)], who discussed Type I error rate, and problems with its assessment, in the evaluation of multiplicity in medical-research literature.

**2. Error rate and false-discovery rate.** Suppose we conduct $\nu$ tests, based on the respective values of the random variables $X_1, \ldots, X_\nu$. Here, $X_i$ typically represents a test statistic computed from the $i$th of a sequence of samples. We reject the $i$th null hypothesis, $H_{0i}$, representing, for example, the hypothesis that the "center" (e.g., the mean) of the population from which the $i$th sample is drawn equals zero, if $X_i > t$; if $X_i \leq t$, then we do not reject $H_{0i}$. Let $N$, a random variable, denote the number of rejected hypotheses:

$$(2.1) \qquad N = \sum_{i=1}^{\nu} I(X_i > t).$$

If each of $H_{01}, \ldots, H_{0\nu}$ is correct, and if we view the sequence of $\nu$ tests as a test of the simultaneous hypothesis $H_0$ that each of the component hypotheses $H_{0i}$ is true, then the significance level of the simultaneous test equals the probability that $N \geq 1$. This is the family-wise error rate (FWER) of the procedure. For example, if $0 < \alpha < 1$ and we define $\beta = -\log(1-\alpha)$; if we choose $t$, in (2.1), to satisfy

$$(2.2) \qquad P_0(X > t) = \frac{\beta}{\nu} + o(\nu^{-1})$$

and if

(2.3)   the random variables $X_i$ are independent and identically distributed as $X$;

then the family-wise error rate converges to $\alpha$ as $\nu$ increases: $P_0(N \geq 1) \to \alpha$. Here and in (2.2), $P_0$ denotes probability computed under $H_0$.

The assumption in (2.3) that the test statistics $X_i$ are identically distributed can be relaxed without much difficulty. For example, if $X_i$ is a Student's $t$-statistic, then it is permissible for $X_i$ to be based on a sample of size $n_i$ drawn from a distribution $F_i$, both depending on $i$, provided the sample sizes and distributions do not vary too greatly with $i$. However, the assumption of independence in (2.3) is critical to our argument at this point.

More generally, it is of interest to determine the probability that we make at least $k$ false discoveries, that is, $P_0(N \geq k)$, where $k \geq 1$ can be arbitrary. This is the generalized family-wise error rate (GFWER). If $t$ satisfies (2.2), and if (2.3) holds, then $N$ is asymptotically Poisson-distributed with mean $\beta$, and so as $\nu \to \infty$,

$$(2.4) \qquad P_0(N \geq k) \to \sum_{j=k}^{\infty} \frac{\beta^j}{j!} e^{-\beta}.$$

An alternative, false-discovery rate (FDR) approach, developed by Simes (1986), Hommel (1988), Hochberg (1988) and Benjamini and Hochberg (1995), involves a step-down procedure but can be framed in a similar way to GFWER. [See also Sarkar (1998) and Sen (1999).] In particular, for $i \geq 1$ let $t_1 > t_2 > \cdots$ denote a sequence depending on $\nu$ and with the property, analogous to (2.2), that

$$P_0(X > t_i) = \frac{i\beta}{\nu} + o(\nu^{-1}). \tag{2.5}$$

[Thus, $t$ in (2.2) is here denoted by $t_1$.] Write $N_i$ for the number of values $X_i$ that lie in the interval $(t_i, t_{i-1}]$, where we take $t_0 = \infty$. The event that the step-down method of Benjamini and Hochberg (1995) does not reject any of the hypotheses $H_{0i}$, for $1 \leq i \leq k$, is equivalent to the event that, for each $i$ in the latter range, $X_i = X_{(\nu-j+1)} \leq t_j$, where $X_{(1)} \leq \cdots \leq X_{(\nu)}$ represent the order statistics of the sequence $X_1, \ldots, X_\nu$. In particular, if $k$ denotes the largest $j$ for which $X_{(\nu-j)} \leq t_{j-1}$, then $H_{0i}$ is rejected for each $i$ such that $X_i = X_{(\nu-j+1)}$, where $1 \leq j \leq k$.

This indicates that, to describe properties of the false-discovery rate approach, we need to understand not just the distribution of $N$, defined at (2.1), but more generally the distribution of

$$N^{(k)} = \sum_{i=1}^{\nu} I(X_i > t_k).$$

Note that $N^{(k)} = N_1 + \cdots + N_k$, where

$$N_i = \sum_{j=1}^{\nu} I(t_i \leq X_j < t_{i-1}). \tag{2.6}$$

Assuming that both (2.3) and (2.5) hold, the random variables $N_1, \ldots, N_k$ are asymptotically independent and Poisson-distributed with mean $\beta$. Therefore, the probability that the null hypotheses corresponding to the $k$ largest values of $X_i$ are all rejected under the FDR approach, when they are in fact all correct, is given by

$$P_0(N^{(i)} \geq i \text{ for } 1 \leq i \leq k) \to P(Q_1 + \cdots + Q_i \geq i \text{ for } 1 \leq i \leq k), \tag{2.7}$$

where $Q_1, \ldots, Q_k$ are independent and identically Poisson-distributed with mean $\beta$. It can be shown from the lemma of Benjamini and Hochberg (1995), page 293, that the probability on the right-hand side of (2.7) is dominated by $\beta$, for each $k \geq 1$. Of course, this is useful only if $\beta < 1$.

In conventional treatments of error rate and false-discovery rate problems, the right-hand sides of (2.2) and (2.5) would generally be replaced by $1 - (1 - \beta)^{1/\nu}$ and $i\beta/\nu$, respectively, reflecting an assumption that the null distribution of $X$ is known exactly. By way of comparison, (2.2) and (2.5) countenance a certain amount of error in our knowledge of the distribution.

The key approximation properties needed to interpret GFWER and FDR in practice are (2.4) and (2.7), which describe the probability of making at least $k$ false discoveries when using the respective methods. In both cases the assumption of independence, in (2.3), is crucial; without it the Poisson approximations may be poor. Our aim is to explore the extent to which the approximations can be rendered invalid by dependence. The context of family-wise error rate is relatively transparent, and so we shall pay greatest attention to that, although giving explicit results in the setting of false-discovery rate.

### 3. Conditions under which clustering occurs, or fails to occur.

3.1. *Models for clustering and for the process $X_i$.* If tests of the hypotheses $H_{0i}$ are conducted independently of one another, then there is no evidence of clustering of level exceedances. In particular, if the random variables $X_i$ are independent and have infinite upper tails, then, trivially,

(3.1)  for each $i_0 \geq 1$     $P(X_i > x$ for some $i$ with $1 \leq |i| \leq i_0 \mid X_0 > x) \to 0$

as $x \to \infty$. We shall define (asymptotic) clustering to occur if (3.1) fails.

Rather than take the $X_i$'s to be independent, we shall model them by a moving average:

$$(3.2) \qquad X_i = \sum_k \theta_k \varepsilon_{i+k},$$

where the $\theta_k$'s are constants and the random variables $\varepsilon_i$, for $-\infty < i < \infty$, are independent and identically distributed. Motivated by simplicity, and by the fact that our definition of clustering involves only fixed, finite values of $i_0$ in (3.1), we shall take the moving average to be of finite order:

$\theta_k = 0$ for all but a finite number of values of $k$, and $\theta_k \neq 0$ for some $k$.
(3.3)

Of course, all our results can be extended to the setting of infinite-order moving averages with sufficiently rapidly decreasing weights $\theta_k$, and in particular, all of the $\theta_k$'s can be nonzero. We confine attention to the finite-order case only for convenience.

The model (3.2) is admittedly rudimentary. However, a more detailed treatment, starting from a "time series" model for the data and, through that, constructing a model for the statistics $X_i$, requires specific information about the definition of the test statistic. The choice at (3.2) is appropriate if the test is being conducted about a mean when the variance is known, and in particular if $X_i = n^{-1/2} \sum_{1 \leq j \leq n} V_{ij}$, where

(3.4)     $V_{ij} = \mu_i + \sum_k \theta_k \varepsilon'_{i+k,j}$     for $1 \leq i \leq \nu$   and   $1 \leq j \leq n$,

$\mu_i = E(V_{ij})$ and the disturbances $\varepsilon'_{ij}$ are all independent and identically distributed with zero expected value. Here, (3.2) holds if we take

$$(3.5) \qquad \varepsilon_i = n^{-1} \sum_{1 \leq j \leq n} \varepsilon'_{ij},$$

these variables being independent and identically distributed. The null and alternative hypotheses under test using the statistic $X_i$ are $H_{0i} : \mu_i = 0$ and $H_{1i} : \mu_i > 0$, respectively.

3.2. *Sufficient conditions for no clustering.* We first state a simple, sufficient condition for (3.1). Let the linear process $X_i$ be as at (3.2), let $\mathcal{K}_j$ denote the set of integers $k$ such that $\theta_{k-j} \neq 0$, and put $\mathcal{K}^{(j)} = \mathcal{K}_j \cap \mathcal{K}_0$. We ask that the independent and identically distributed disturbances $\varepsilon_i$ satisfy

$$(3.6) \text{ for each } v > 0 \text{ and each } j \neq 0 \qquad \frac{P(\sum_{k \in \mathcal{K}^{(j)}} \theta_k \varepsilon_k > u - v)}{P(\sum_{k \in \mathcal{K}_0} \theta_k \varepsilon_k > u)} \to 0$$

$$\text{as } u \to \infty.$$

Let $1 \leq I_1 < I_2 < \cdots$ denote the indices $i$ for which $X_i > t$, where $t$ is as in (2.2).

THEOREM 3.1. *If (3.3) and (3.6) hold, then so too does (3.1).*

THEOREM 3.2. *If (2.2) and (3.1) hold, then, for each constant $C > 0$, the point process $I_1\nu^{-1}, I_2\nu^{-1}, \ldots$, restricted to the interval $[0, C]$, converges weakly, as $\nu \to \infty$, to a homogeneous Poisson process on $[0, C]$, with intensity $\beta$.*

Theorem 3.2 implies that $N$, at (2.1), is Poisson-distributed with mean $\beta$. The argument leading to Theorem 3.2 also shows that, if (2.5) and (3.1) hold, then for each $i \geq 1$ the random variables $N_1, \ldots, N_i$, introduced at (2.6), are independent and identically Poisson-distributed with mean $\beta$. Together these results establish the correctness of the crucial Poisson approximations (2.4) and (2.7).

As noted in Section 2, these results also hold if $X_1, X_2, \ldots$ are independently distributed. Therefore, under conditions (2.5) and (3.1), exceedences of the level $t$ by the linear process $X_1, X_2, \ldots$ have the same first-order asymptotic properties they would enjoy if the $X_i$'s were independent and identically distributed random variables with the same marginal distribution as the linear process. In particular, the Introduction of dependence does not produce any first-order evidence of clustering.

Therefore, calibrating the tests using methodology based on the assumption of independence is adequate if the null distribution of the stochastic process $X_i$ is close to that of a linear process, if the number of simultaneous tests is sufficiently large, and if (3.1) holds. In the next section we shall show that (3.6), and hence (3.1), prevails if the marginal distribution of $X_i$ is light-tailed.

3.3. *No clustering occurs for light-tailed distributions.* Here we show that, under the moving-average model defined at (3.2) and (3.3), no clustering occurs [i.e., (3.1) holds] if the distribution tails decrease like $\exp(-\operatorname{const.} x^\gamma)$ where $\gamma \geq 1$. Therefore, testing can proceed as though the test statistics $X_i$ are independent, which of course they are not.

The case where $\gamma > 1$ is relatively straightforward; there we need assume only that, for a constant $C > 0$, the density $f$ of the distribution of $\varepsilon$ satisfies, as $x \to \infty$,

$$(3.7) \qquad f(x) = \exp\{o(x^\gamma)\} \exp(-Cx^\gamma).$$

A sufficient condition for (3.7) is the following: For constants $C, C_1 > 0$ and $C_2 \geq 0$,

$$(3.8) \qquad f(x) \sim C_1 x^{C_2} \exp(-Cx^\gamma)$$

as $x \to \infty$.

THEOREM 3.3. *If the process $X_1, X_2, \ldots$ is determined by (3.2); if the density of $\varepsilon$ exists and satisfies (3.7) with $\gamma > 1$, or satisfies (3.8) with $\gamma = 1$; and if the weights $\theta_k$ are all nonnegative and satisfy (3.3); then (3.6), and hence also (3.1), hold.*

The assumption that the weights $\theta_k$ are all nonnegative is important, in that without it, properties of the lower tail of the distribution of $\varepsilon$ would have to be taken into account. [Conditions (3.7) and (3.8) address only the upper tail.] Depending on behavior of the lower tail, if one or more of the $\theta_k$'s is negative, then first-order asymptotic theory can be quite different from that discussed in Theorems 3.3–3.5.

For example, if the negative $\theta_k$'s form a set $\{\theta_k \equiv -\omega, \text{ for } k \in \mathcal{A}\}$, where $\omega > 0$; and if the density of the lower tail of the distribution of $\varepsilon$ satisfies

$$f(-x) \sim C_3 x^{C_4} \exp(-C_5 x^{\gamma_1})$$

as $x \to \infty$, where $C_3, C_5 > 0$, $C_4 \geq 0$ and $0 < \gamma_1 < 1 < \gamma$; then the pattern of exceedences of $t$ [where $t$ still has the property at (2.2)] is first-order equivalent to that for quite a different process $X_i$, for which the only nonzero moving-average weights are $\theta_k \equiv \omega$ for $k \in \mathcal{A}$, and where the distribution of $\varepsilon$ satisfies (3.8) with $(C, C_1, C_2, \gamma)$ there replaced by $(C_5, C_3, C_4, \gamma_1)$. For such a process, clustering can occur; see Theorem 3.4 below. Thus, by allowing negative weights and choosing the lower-tail distribution appropriately, we can substantially alter the pattern of level exceedances.

The case of Student's $t$-statistic is related to the model (3.2), but differs in important respects. One of these is the potential for the tails of the distribution of $X_i$ to become lighter as the "group size," that is, the size of the dataset used to compute an individual $X_i$, increases. We shall discuss this issue in Section 3.6.

Theorem 3.3 includes the case where the autoregression is a Gaussian process. In particular it implies that, in the Gaussian setting, clustering does not occur unless, for example, the strength of dependence of the process $X_i$ is permitted to increase with $\nu$. We shall take up this issue in Section 3.7, showing that correlations must converge to 1 at least as fast as $(\log \nu)^{-1}$ if clustering is to be present in asymptotic terms.

3.4. *Clustering can sometimes occur if* $0 < \gamma < 1$. The case where (3.7) or (3.8) holds, and $0 < \gamma < 1$, is relatively complex. There, if the largest $\theta_k$ occurs for a unique value of $k$, then (3.1) holds. That is, the probability that there exists a cluster of exceedances converges to zero as the exceedance level, $x$, increases. In this instance, if $\nu$ is sufficiently large, the dependent test statistics $X_i$ can be treated as though they were independent, without serious problems arising.

However, if there are ties for the largest $\theta_k$, then the probability of a cluster does not converge to zero. In this case, if the number of tied values equals $q$, then the probability that the size of the cluster of exceedances also equals $q$, converges to 1 as the exceedance level increases.

To indicate why the case $\gamma < 1$ is so different, we treat the instance where $\theta_1 = \cdots = \theta_r$ and each other $\theta_k$ vanishes. In this setting, having $\varepsilon_1 + \cdots + \varepsilon_r > x$ implies that, with high probability, one of the values of $\varepsilon_1, \ldots, \varepsilon_r$ is very close to $x$, or greater than $x$, and the other values are all significantly smaller than $x$. (Here and below we assume that $x$ is large.) That is, just one of the $\varepsilon_i$'s is responsible for the level exceedance, and its influence can persist, through weights in the moving average, to ensure that $\varepsilon_{j+1} + \cdots + \varepsilon_{j+r} > x$ for values of $j$ other than simply $j = 0$.

By way of comparison, if $\gamma > 1$ and $\varepsilon_1 + \cdots + \varepsilon_r > x$, then it is highly likely that this is achieved through all of the $\varepsilon_k$'s being of order $x$; and that, if just one of the $\varepsilon_i$'s is exchanged for another, the inequality fails. Therefore in this case, if $\varepsilon_1 + \cdots + \varepsilon_r > x$, then it is unlikely that $\varepsilon_{j+1} + \cdots + \varepsilon_{j+r} > x$ for values of $j \neq 0$. Therefore clustering can occur if $\gamma < 1$, but is relatively unlikely if $\gamma > 1$. Our proofs in Section 5 involve verification of general versions of these properties, which underpin the intuitive arguments given in Section 1.

Next we formally state a result describing the case $0 < \gamma < 1$. Write $r$ for any integer that is not less than the difference between the least, and largest, values of $k$ for which $\theta_k \neq 0$, and let $M$ denote the number of values $j$ with $|j| \leq r$, for which $X_j > x$.

THEOREM 3.4. *Assume that* (a) *the weights $\theta_k$ are all nonnegative and satisfy (3.3), and* (b) *the density $f$ of the distribution of $\varepsilon$ exists and satisfies (3.7) for a value of $\gamma$ in the range $0 < \gamma < 1$. If, in addition,* (c) *there is no tie for the largest $\theta_k$, then* (i) *(3.1) holds. On the other hand, if* (a) *and* (b) *hold, although with (3.8) replacing (3.7) in* (b) *and, instead of* (c), (d) *exactly $q \geq 2$ of the values of $\theta_k$ tie for the maximum, then* (ii) *$P(M = q \mid X_0 > x) \to 1$ as $x \to \infty$.*

It follows from Theorems 3.2 and 3.4 that if (2.5) and (a)–(c) in Theorem 3.4 hold, then the random variable $N$, at (2.1), is asymptotically Poisson with mean $\beta$; and likewise, that the random variables $N_1, \ldots, N_k$, defined at (2.6), are asymptotically independent and Poisson with mean $\beta$. This shows that, asymptotically, clusters do not occur, and establishes the correctness of the key Poisson approximations, (2.4) and (2.7), borrowed from the case where the $X_i$'s are independent.

However, if (c) in Theorem 3.4 fails, and is replaced there by (d), then with probability converging to 1, clusters exist and are of size $q$. Moreover, $q^{-1}N$ is asymptotically Poisson, and $q^{-1}N_1, \ldots, q^{-1}N_k$ are asymptotically independent and Poisson, with mean $\beta/q$ in each case. Therefore (2.4) and (2.7) fail in this case. For example, (2.7) should be replaced by the result,

$$P_0(N^{(i)} \geq i \text{ for } 1 \leq i \leq k) \to P(qQ_1 + \cdots + qQ_i \geq i \text{ for } 1 \leq i \leq k),$$

where $Q_1, \ldots, Q_k$ are independent and Poisson with mean $\beta/q$. The fact that $q^{-1}N$ and $q^{-1}N_i$, rather than $N$ and $N_i$, are independent and Poisson, follows using part (ii) of Theorem 3.4 and the fact that the probability that a cluster overlaps the end of the interval $1, 2, \ldots, \nu$ converges to zero as $\nu \to \infty$.

To appreciate intuitively why, in the paragraph above, the Poisson mean equals $\beta/q$ rather than $\beta$, note that (2.2) and (2.5) imply that $\nu E(N) \to \beta$ and $\nu E(N_i) \to \beta$ as $\nu \to \infty$. However, each time an exceedence occurs it is, with probability converging to 1, accompanied by $q - 1$ other exceedences, and so if the number of clusters has mean $\beta_1$, then $\beta_1 q = \beta$, that is, $\beta_1 = \beta/q$.

3.5. *Clustering in the case of Pareto-type distributions of disturbances.* Here we assume that, for constants $C, \rho > 0$,

$$(3.9) \qquad P(\varepsilon > x) \sim C x^{-\rho}$$

as $x \to \infty$. More generally, $C$ could be replaced by a slowly varying function of $X$. In these settings the probability that a cluster of exceedences occurs is bounded away from zero, as the exceedence level increases, regardless of ties among the moving-average weights.

To describe the distribution of cluster size, let $\theta_{(1)} \geq \cdots \geq \theta_{(m)}$ denote a ranking of the $m$ nonzero $\theta_i$'s. Define $\theta_{(q)} = 0$ for $q > m$ and $p_q = (\theta_{(q)}^\rho - $ $\theta_{(q+1)}^{\rho})/\theta_{(1)}^{\rho}$. Let $M_0$ denote a random variable for which $P(M_0 = q) = p_q$. Note that if all the nonzero $\theta_i$'s are equal, then $P(M_0 = m) = 1$. Our next theorem asserts that the distribution of $M_0$ is the limiting distribution of cluster size. Given $x > 0$, write $M$ for the number of values $j$ with $|j| \leq r$, for which $X_j > x$, and define $M_1$ to have the distribution of $M$ given that $M \geq 1$.

THEOREM 3.5.    *If (3.9) holds, and if the weights $\theta_k$ are all nonnegative and satisfy (3.3), then $P(M_1 = q) \to P(M_0 = q)$ as $x \to \infty$.*

Theorem 3.5 implies that both (2.4) and (2.7) fail in the forms given there. We now outline modifications to (2.4) and (2.7) that are necessary if those results are to hold in the setting of (3.9).

Put $\mu = E(M_0)$, let $Q$ and $Q_1, Q_2, \ldots$ be independent and identically Poisson-distributed random variables with mean $\beta/\mu$, and let $M_1, M_2, \ldots$ and $M_{j\ell}$, for $j \geq 1$ and $\ell \geq 1$, be independent random variables each with the distribution of $M_0$. In cases where (3.9) holds, (2.4) and (2.7) should be replaced by, respectively,

$$(3.10) \qquad P_0(N \geq k) \to P\left(\sum_{i=1}^{Q} M_i \geq k\right),$$

$$(3.11) \quad P_0(N^{(i)} \geq i \text{ for } 1 \leq i \leq k) \to P\left(\sum_{j=1}^{i} \sum_{\ell=1}^{Q_j} M_{j\ell} \geq i \text{ for } 1 \leq i \leq k\right).$$

In principle the Pareto parameter $\rho$, and the constants $\theta_k$ in the linear-process model, also can be estimated from data, and hence the distribution of $M_0$ can be estimated. This leads to estimators of the right-hand sides of (3.10) and (3.11). However, this approach to statistical analysis will generally not be straightforward.

3.6. *The case of Student's t-statistic.*    The model (3.2) for $X_i$ is directly appropriate when the test statistic is a sample mean, but in other cases it is only an approximation. For example, in the context of two-channel microarrays, $X_i$ would be a Studentized mean. In this setting, suppose data $V_{i1}, \ldots, V_{in}$ are generated as at (3.4), and consider the test that rejects $H_{0i} : \mu_i = 0$, in favor of $H_{1i} : \mu_i > 0$, if $Y_i > t$, where

$$(3.12) \qquad Y_i = \frac{n^{-1/2} \sum_{1 \leq j \leq n} V_{ij}}{\{n^{-1} \sum_{1 \leq j \leq n} V_{ij}^2 - (n^{-1} \sum_{1 \leq j \leq n} V_{ij})^2\}^{1/2}}$$

is a conventional $t$-statistic. If $n$ is large, then the distribution of $Y_i$ under $H_{0i}$ can be approximated by the distribution of $X_i$, at (3.2), on taking $\varepsilon_i$ to be given by (3.5). Moreover, as $n$ increases the distribution of $Y_i$ becomes more light-tailed, and so high-level exceedences by the $Y_i$'s should become less clustered. Perhaps surprisingly, "large" $n$ can be very much less than $\nu$ [it is sufficient that $\log \nu = o(n)$ as $n$ diverges], and the tails of the distribution of $\varepsilon$ can be relatively heavy (only $E|\varepsilon|^3 < \infty$ is required), without damaging the property that high-level crossings are asymptotically independent. Also, depending on the weights $\theta_k$, the level, $t$, at which these properties occur can be substantially lower than in the setting of Theorems 3.1–3.3. These results make substantial use of special properties of $t$-statistics, and will be given elsewhere.

3.7. *The case of a highly correlated Gaussian process.* The reader will have noticed that the strength of dependence permitted by the model (3.2) is reasonably low, and might well ask: "Just how strong does dependence have to be before clustering becomes apparent?" Our purpose in Section 3.7 is to respond to that question. In the context of processes for which dependence decays to zero over a finite range, the answer is, "The point at which clustering is noticed is where the correlation between nearby $X_i$'s is $1 - \text{const.}(\log \nu)^{-1/2} + o\{(\log \nu)^{-1/2}\}$." This is not especially strong correlation; for each $\eta > 0$ it is weaker than $1 - \text{const.} \nu^{-\eta}$.

There exist real-world processes where dependence at neighboring indices $i$ can be very strong. Consider, for example, the case of speckle imaging in astronomy, where noise correlation at neighboring pixels can be particularly high. This has a significant effect on the potential for resolving (or for successfully testing for the existence of) faint light sources in the heavens.

To model these processes we shall take the variables $\varepsilon_i$, in the moving average at (3.2), to be independent $N(0,1)$ random variables, and the weights $\theta_k$ to be given by

$$(3.13) \qquad \theta_{-k} = c \prod_{j=0}^{k} \rho_k \qquad \text{for } k \geq 0, \theta_k = 0 \qquad \text{for } k \geq 1,$$

where the constants $\rho_k$ are nonnegative, and $c > 0$ is chosen so that var $X_i = 1$ for each $i$. If each $\rho_k = \rho$, not depending on $k$, then $X_i$ is an autoregression of order 1: $X_i = \rho X_i + (1 - \rho^2)^{1/2} \varepsilon_i$. We shall instead take

$$(3.14) \qquad \begin{aligned} &\rho_0 = 1, \qquad \rho_k = 1 - a_k \delta + o(\delta) \qquad \text{for } 1 \leq k \leq r, \\ &\rho_k = 0 \qquad \text{for } k \geq r + 1, \end{aligned}$$

where $\delta = \delta(\nu) \downarrow 0$ as $\nu \to \infty$, and $a_1, \ldots, a_r$ are nonnegative constants.

Define

$$(3.15) \qquad c_j = \frac{1}{r+1} \sum_{k=0}^{r} (a_{k+1} + \cdots + a_{k+j}).$$

Then, $\text{cov}(X_i, X_{i-j}) = 0$ for $j \geq r+1$, whereas for $0 \leq j \leq r$,

$$(3.16) \qquad \text{cov}(X_i, X_{i-j}) = 1 - c_j \delta + o(\delta).$$

These properties, and the fact that $\delta$ decreases with increasing $\nu$, imply that $X_{i_1}$ and $X_{i_2}$ are very highly correlated if $|i_1 - i_2| \leq r$, but are independent otherwise.

We shall give the limiting distribution of cluster size in this setting. To do so, define $\mathcal{I}$ to be the set of $2r$ integers between $-r$ and $r$, excluding zero; and let $Z_i$, for $i \in \mathcal{I}$, denote $2r$ Normally distributed random variables with zero means and covariance matrix $\Sigma = (\sigma_{ij})$, where

$$\sigma_{ij} = \text{cov}(Z_i, Z_j) = c_{|i|} + c_{|j|} - c_{|i-j|}$$

and $c_j$ is as at (3.15). Write $\bowtie$ to denote either ">" or "<," and let $\mathcal{S} = (\bowtie_i : i \in \mathcal{I})$ be a sequence of such inequalities. Of course, there are just $2^{2r}$ distinct sequences $\mathcal{S}$. Given a constant $d > 0$, and given a particular sequence $\mathcal{S}$, define

$$\pi(\mathcal{S}) = \int_0^\infty P(Z_i \bowtie_i d c_{|i|} - d^{-1} z \text{ for } 1 \leq |i| \leq r) e^{-z} \, dz.$$

For $0 \leq k \leq 2r$, let $\pi_k^0$ equal the sum of $\pi(\mathcal{S})$ over all sequences $\mathcal{S}$ that contain just $k$ ">" signs and $2r - k$ "<" signs. Define $\pi_k = \pi_k(\nu)$ to equal the probability that exactly $k$ out of the $2r$ values of $X_i$, for $i \in \mathcal{I}$, exceed $t$, conditional on $X_0 > t$.

THEOREM 3.6.   *If the errors $\varepsilon_i$ are independent Normal $N(0,1)$, so that the process $X_i$, defined at (3.2), is Gaussian; if the weights $\theta_k$ are given by (3.13), and the coefficients $\rho_k$ are given by (3.14), with $a_1, \ldots, a_r \geq 0$; if $c_1, \ldots, c_k$ are defined in terms of $a_1, \ldots, a_k$ by (3.15); and if $t$ and $\delta^{-1}$ both diverge as $\nu \to \infty$, with $\delta^{1/2} t \to d$, where $0 \leq d \leq \infty$; then, for $0 \leq k \leq 2r$, $\pi_k \to \pi_k^0$ if $0 < d < \infty$, $\pi_k \to 0$ if $d = \infty$, and $\pi_k \to 1$ if $d = 0$.*

Note that, when $X$ has a normal $N(0,1)$ distribution, the value of $t$ defined by (2.2) satisfies $t \sim (2 \log \nu)^{1/2}$ as $\nu$ increases. Therefore the condition invoked in Theorem 3.6, that $\delta^{1/2} t \to d$ for some finite and nonzero $d$, is equivalent to the correlation between neighboring $X_i$'s equalling $1 - \text{const.} (\log \nu)^{-1/2} + o\{(\log \nu)^{-1/2}\}$.

## 4. Numerical properties.

Our simulations were based on two different models. In model 1 the test statistic $X_i$ was that given at (3.2), with $\varepsilon$ simulated from a Student's $t$ distribution. Model 2 was the Student's $t$-statistic model at (3.12) with $n = 10$; we took the distribution of $\varepsilon$ to itself be Student's $t$. In both models the number of nonzero $\theta_k$'s (which we shall call $r$) was taken to equal 1 (independence), 3, 10 or 50, and the nonzero $\theta_k$'s were taken equal to one another.

The number, $\nu$, of tests was 500, 1000, 2000, 5000 or 10,000 for both models. A range of tail weights was achieved by varying the number of degrees of freedom for the distribution of $\varepsilon$; we included infinity, thereby addressing the case of normally distributed $\varepsilon$. These were scaled so that $\text{var}(X_i) = 1$ in each case. The chosen critical values were based on controlling the FWER in the one-sided case with $\alpha = 0.05$. Each simulation involved 10,000 repetitions.

"Clustering tendency" can be characterized in terms of the value of $N$, that is, the number of rejected hypotheses. If the hypothesis tests are genuinely independent, then most realizations have $N$ equal to 0 or 1; the proportion of realizations for which $N > 1$ is only 0.0013. However, as the effects of dependence become more pronounced, leading to greater clustering, the event $N > 1$ becomes more common, with a corresponding decrease in the number of events for which $N = 1$. Therefore a succinct way of reporting the effect that tail-weight of the error distribution has on clustering tendency is to graph the proportion of clusters for which $N > 1$ of those for which $N > 0$, against number of degrees of freedom (df).

This is the approach taken in Figures 1 and 2, which summarize these results. In both figures, panels (a) through (d) represent the different values of $r$ (1, 3, 10 or 50, resp.). The horizontal axis gives the number of degrees of freedom, and each separate line represents a different number of tests, $\nu$.

As Figure 1 indicates, in the case of model 1 there is a clear decrease in clustering as tail-weight decreases for $r = 3$, 10 and 50. This reflects the results in Theorem 3.3, for example. There is also a slightly less clear, but nevertheless present, decrease in clustering as $\nu$ increases, particularly for normally distributed $\varepsilon$. While these trends are present for all values of $r$, by the time $r$ is as large as 10 the strength of dependence has increased so much that the decrease in clustering with decreasing tail-weight is noticeably slower. See, for example, the panels of Figure 1 corresponding to $r = 10, 50$.

Reflecting the conclusions reached in Section 3.6, Figure 2 indicates that there is very little clustering under model 2 for $r = 3$, even for heavy-tailed $\varepsilon$. There is still clustering for long-range dependency, which persists in the light-tailed case, although it decreases as $\nu$ increases.

The case of nonequal $\theta_k$ was considered; the cases with larger $r$ behaved like those with smaller $r$ if the number of large $\theta_k$'s was small.

## 5. Technical arguments.

5.1. *Proof of Theorem 3.1.* To derive (3.1) it suffices to show that, for each $j$ for which $\mathcal{K}_j$ is a proper subset of $\mathcal{K}_0$, $P(X_j > x \mid X_0 > x) \to 0$. To

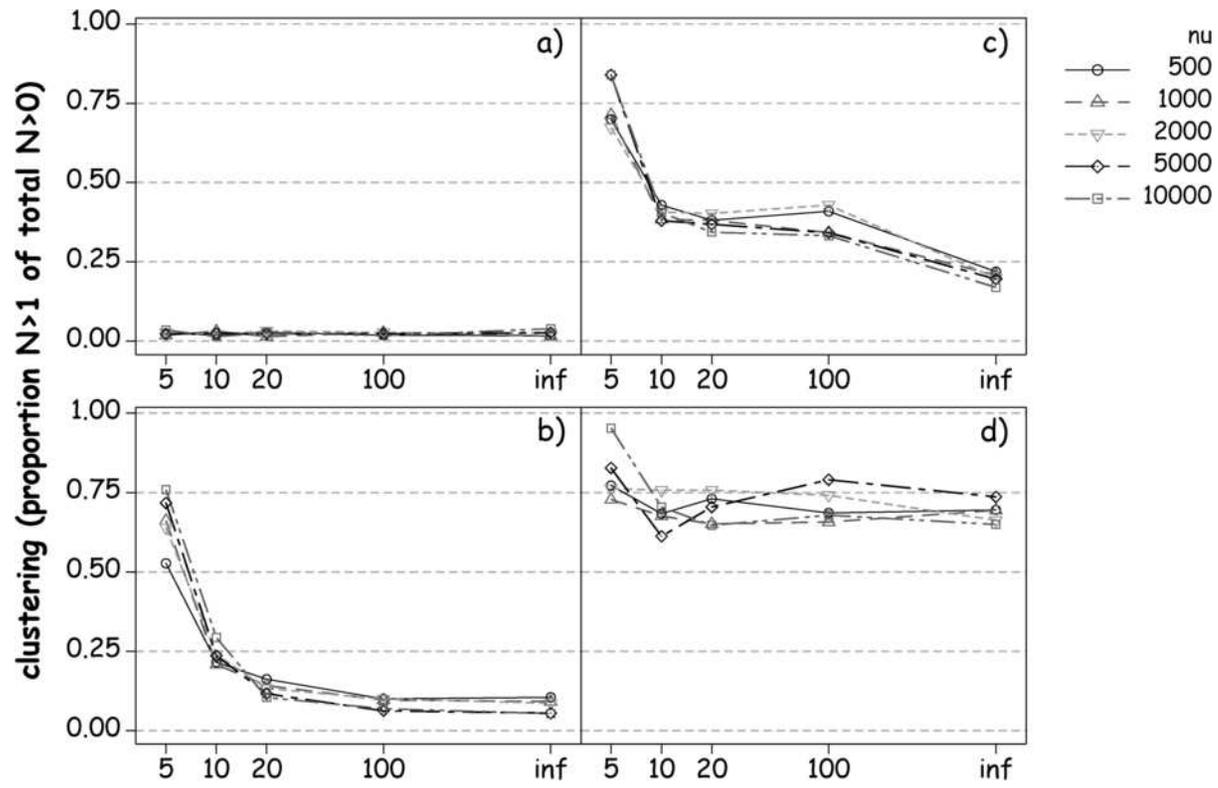

FIG. 1. *Clustering when test statistics are distributed as moving averages. (a) One nonzero value of $\theta_k$; (b) three nonzero values of $\theta_k$; (c) ten nonzero values of $\theta_k$; (d) fifty nonzero values of $\theta_k$ (for clarity, the horizontal axis is logarithmic).*

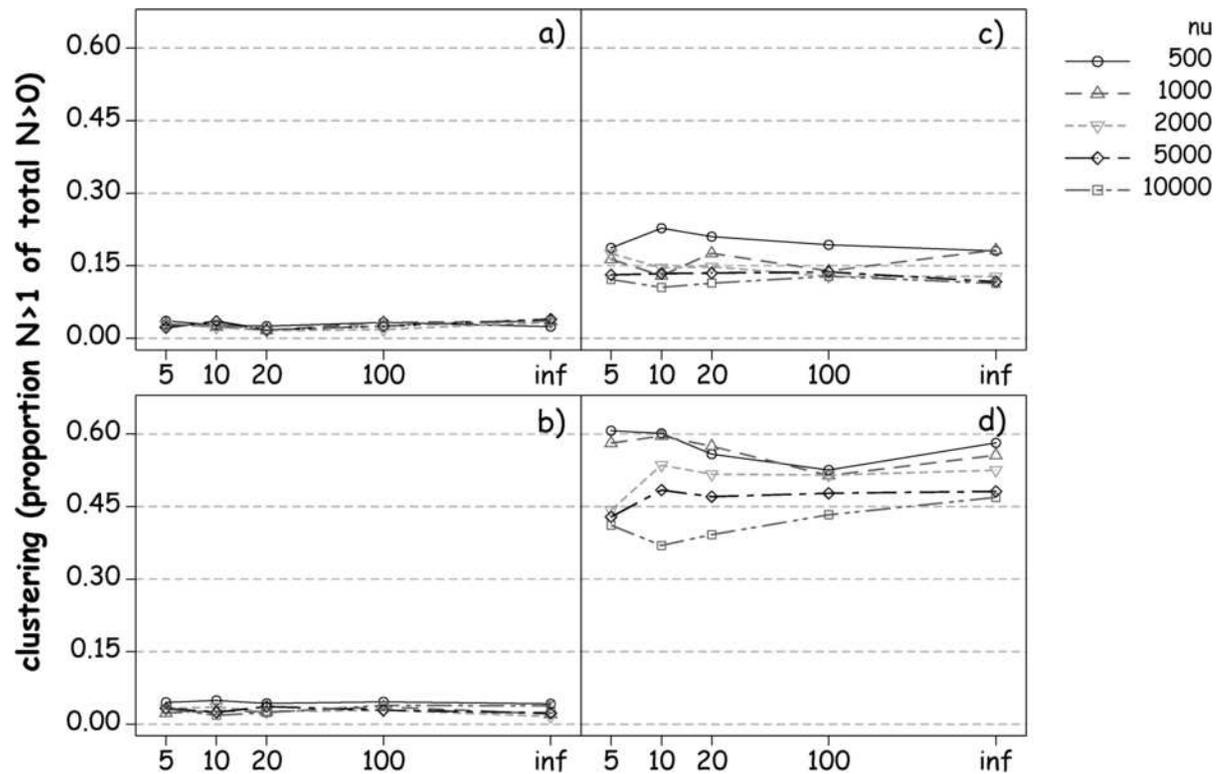

FIG. 2. *Clustering when test statistics are t-statistics computed from moving-average data. (a) One nonzero value of $\theta_k$; (b) three nonzero values of $\theta_k$; (c) ten nonzero values of $\theta_k$; (d) fifty nonzero values of $\theta_k$ (for clarity, the horizontal axis is logarithmic).*

this end, put

$$U = \sum_{k \in \mathcal{K}_j \cap \widetilde{\mathcal{K}}_0} \theta_{k-j} \varepsilon_k, \qquad V = \sum_{k \in \mathcal{K}_j \cap \mathcal{K}_0} \theta_{k-j} \varepsilon_k,$$

where $\widetilde{\mathcal{K}}_0$ denotes the complement of $\mathcal{K}_0$. Then $U$ is independent of both $V$ and $X_0$, and so

(5.1)
$$\begin{aligned}
P(X_j > x \mid X_0 > x) \\
= P(U + V > x \mid X_0 > x) \\
\leq P(u + V > x \mid X_0 > x) + P(U > u) \\
\leq \frac{P(u + V > x)}{P(X_0 > x)} + P(U > u).
\end{aligned}$$

The ratio $P(u + V > x)/P(X_0 > x)$ has the same form as the ratio of probabilities in (3.6), with $(u, v)$ there replaced here by $(x, u)$. Hence, by (3.6), the far right-hand side of (5.1) converges to $P(U > u)$ as $x \to \infty$. Since this is true for arbitrarily large $u$, then the far left-hand side of (5.1) converges to zero as $x \to \infty$. This proves (3.1).

5.2. *Proof of Theorem 3.2.* Let the integer $\ell$ be so large that, for some $j$, the only $\theta_k$'s for which $\theta_k \neq 0$ are included in the set $\theta_{j+1}, \ldots, \theta_{j+\ell}$; and let $m = m(\nu) \geq 1$ denote an integer satisfying $m \asymp \nu$ as $\nu \to \infty$. Divide the indices $1, \ldots, m$ into $B$ blocks, each of length $b = b(\nu)$, where $b \to \infty$ and $b/\nu \to 0$ as $\nu$ increases; with consecutive blocks separated by "spacers" of length $\ell$; in such a way that $X_1, \ldots, X_b$ and $X_{m-b+1}, \ldots, X_m$ denote the first and last block, respectively. (This neat fit of the blocks, and their separating spacers, into the interval $[1, m]$ may require a slight increase in $m$, but since $b/\nu \to 0$, then the fit may be achieved without damaging the property $m \asymp \nu$.) Define $J_j = 1$ (resp., $K_j = 1$) if $X_i > t$ for some integer $i$ in the $j$th block (in the $j$th spacer), and put $J_j = 0$ ($K_j = 0$) otherwise. Then $J_1, J_2, \ldots$ are independent random variables [call this property (P$_1$)], as too are $K_1, K_2, \ldots$. Let $\mathcal{J}(b)$ denote the set of indices $j$ such that $\mathcal{K}^{(j)}$ is a proper subset of $\mathcal{K}_0$ and $|j| \leq b$. If $b$ diverges to infinity sufficiently slowly, then, by (3.1), and as $\nu \to \infty$,

(5.2)        $P\{X_j > t \text{ for some } j \in \mathcal{J}(b) \mid X_0 > t\} \to 0.$

Let $M_J$ and $M_K$ denote the number of nonzero $J_j$'s, and number of nonzero $K_j$'s, respectively. Markov's inequality, (2.2), and the fact that $b \to \infty$ as $\nu \to \infty$, can be used to show that as $\nu$ increases, $P(M_K = 0) \to 1$ [call this property (P$_2$)] and

(5.3)        $$\lim_{k \to \infty} \limsup_{\nu \to \infty} P(M_J \geq k) = 0.$$

Result (5.3) implies that the probability that $M_J \leq k$ can be made arbitrarily close to 1, uniformly in $\nu$, by choosing $k$ sufficiently large but fixed. This property, and (5.2), imply that, with probability converging to 1 as $\nu \to \infty$, none of the blocks enjoys more than a single exceedence; call this property ($P_3$). Together, ($P_2$) and ($P_3$) imply that, with probability converging to 1 as $\nu \to \infty$, the number of indices $i$, for $1 \leq i \leq m$, such that $X_i > t$, equals the number of indices $j$, for $1 \leq j \leq B$, such that $J_j = 1$. Call this property ($P_4$).

The Poisson property stated in the theorem, but for the interval $[0, m/\nu]$ rather than $[0, C]$, follows from ($P_1$) and ($P_4$). By taking $m = m(\nu)$ to be so large that $m/\nu \geq C$ for all sufficiently large $\nu$, we complete the proof of the theorem. This argument does not immediately give the mean of the Poisson distribution. However, simple calculations from (2.2) show that $P(J_j = 1) = b\beta\nu^{-1} + o(b\nu^{-1})$, not depending on $j$. From this result, and the fact that $B \sim \nu/b$, follows the claim in the theorem that the limiting Poisson process has intensity $\beta$.

5.3. *Proof of Theorem 3.3.* First we assume that $\gamma > 1$. Without loss of generality, the constant $C$ in (5.3) equals 1. We shall prove that in this case, if $\theta_1, \ldots, \theta_r$ are nonnegative constants, at least one of them positive, and if $\varepsilon_1, \ldots, \varepsilon_r$ are independent and identically distributed random variables for which the density satisfies (3.7), then

$$P(x) \equiv P\left(\sum_{k=1}^{r} \theta_k \varepsilon_k > x\right) = \exp\left\{-\left(\sum_{k=1}^{r} \theta_k^{\gamma/(\gamma-1)}\right)^{-(\gamma-1)} x^\gamma + o(x^\gamma)\right\}.$$

(5.4)

Result (3.6) follows directly.

Let $\mathcal{U}$ denote the set of points $(u_1, \ldots, u_r)$ such that $\sum_k \theta_k u_k > 1$ and each $u_k \geq 0$. It can be deduced from (3.7) that, as $x \to \infty$,

$$P(x) = \exp\{o(x^\gamma)\} \int_{\mathcal{U}} \exp\{-(u_1^\gamma + \cdots + u_r^\gamma)x^\gamma\} \, du_1 \cdots du_r.$$

A Lagrange multiplier argument shows that the minimum of $u_1^\gamma + \cdots + u_r^\gamma$, subject to $\sum_k \theta_k u_k \geq 1$ and each $u_k \geq 0$, occurs when $u_k = C_1 \theta_k^{1/(\gamma-1)}$ and equals $C_1^{\gamma-1}$, where $C_1^{-1} = \sum_k \theta_k^{\gamma/(\gamma-1)}$. Therefore, (5.4) holds.

Next we treat the case $\gamma = 1$. It suffices to assume that $C = 1$ and $P(\varepsilon > 0) = 1$. Suppose too that, among the positive weights $\theta_1, \ldots, \theta_r$, there are just $d$ distinct values of $\theta_k$, given by $\omega_1 > \cdots > \omega_d > 0$, and that these are repeated $s_1, \ldots, s_d$ times, respectively. Thus, $s_1 + \cdots + s_d$ equals the number, $r$, of integers $k$ for which $\theta_k$ is nonzero. Then, writing $du$ for either $du_1 \cdots du_r$ or $du_1 \cdots du_d$, depending on occasion, we have

$$p_1(x) \equiv P\left(\sum_{k=1}^{r} \theta_k \varepsilon_k > x\right)$$

$$\asymp \int_{\substack{\theta_1 u_1 + \cdots + \theta_r u_r > x, \\ u_1, \ldots, u_r > 0}} \left( \prod_{k=1}^{r} u_k^{C_2} \right) \exp\{-(u_1 + \cdots + u_r)\} \, du$$

$$\asymp \int_{\substack{\omega_1 u_1 + \cdots + \omega_d u_d > x, \\ u_1, \ldots, u_d > 0}} \left( \prod_{k=1}^{d} u_k^{s_k(C_2+1)-1} \right) \exp\{-(u_1 + \cdots + u_d)\} \, du$$

$$= \int_{\substack{u_1 > (x/\omega_1) - (\omega_2 u_2 + \cdots + \omega_r u_r)/\omega_1, \\ u_1, \ldots, u_d > 0}} \left( \prod_{k=1}^{d} u_k^{s_k(C_2+1)-1} \right)$$
$$\times \exp\{-(u_1 + \cdots + u_d)\} \, du$$

(5.5)

$$\asymp x^{s_1(C_2+1)-1}$$
$$\times \int_{\substack{\omega_2 u_2 + \cdots + \omega_r u_r \leq x, \\ u_2, \ldots, u_d > 0}} \left( \prod_{k=2}^{d} u_k^{s_k(C_2+1)-1} \right)$$
$$\times \exp[-\{x\omega_1^{-1} + (1 - \omega_1^{-1}\omega_2)u_2$$
$$+ \cdots + (1 - \omega_1^{-1}\omega_d)u_d\}] \, du_2 \cdots du_d$$

$$\asymp x^{s_1(C_2+1)-1} \exp(-x/\omega_1).$$

Result (5.5) gives an asymptotic expression for the denominator in (3.6). An asymptotic formula for the probability in the numerator, equal to $p_2(u - v)$ say, can be derived similarly. To appreciate the conclusion of those calculations, let $\mathcal{K}(j)$ denote the set of indices $k$ in $\mathcal{K}_j$ for which $\theta_k = \omega_1$, and put $\mathcal{K}[j] = \mathcal{K}(j) \cap \mathcal{K}(0)$. Then, for $j \neq 0$, $\mathcal{K}[j]$ is a proper subset of $\mathcal{K}[0] = \mathcal{K}(0)$. If $\mathcal{K}[j]$ is empty, then $p_2(x) = O(e^{-x/\omega})$ for a constant $\omega \in (0, \omega_1)$. If $\mathcal{K}[j]$ contains at most $s_1 - 1$ ($\geq 1$) elements, then $p_2(x) = O\{x^{s_1(C_2+1)-2} \exp(-x/\omega_1)\}$. It follows from these properties and (5.5) that, for each $v$, $p_2(u - v)/p_1(u) \to 0$ as $u \to \infty$. Therefore (3.6) holds. This establishes Theorem 3.3 in the case $\gamma = 1$.

5.4. *Proof of Theorem 3.4.* In order to prove (i) it suffices to show that, for each $j \neq 0$,

(5.6)
$$\frac{P(X_0 > x, X_j > x)}{P(X_0 > x)} \to 0$$

as $x \to \infty$. To achieve this end we shall derive upper and lower bounds for the numerator and denominator, respectively, on the left-hand side.

Let $r$ denote the number of nonzero values of $\theta_k$, choose $B > 0$ so large that $a \equiv P(0 \leq \varepsilon \leq B) > 0$, write $\omega_1$ and $\omega_2$ for the largest and second-largest, respectively, values of $\theta_k$, and put $C_3 = (r-1)\omega_2 B/\omega_1$. Then, since

$0 < \gamma < 1$,

$$P(X_0 > x) \geq P\{\omega_1 \varepsilon > x - (r-1)\omega_2 B\}a^{r-1}$$

(5.7)
$$= \exp\{-C(x\omega_1^{-1} - C_3)^\gamma + o(x^\gamma)\}$$

$$= \exp\{-C(x/\omega_1)^\gamma + o(x^\gamma)\}.$$

Next we derive an upper bound to the numerator in (5.6). We may assume, without loss of generality, that $X_i = \theta_1 \varepsilon_{i+1} + \cdots + \theta_r \varepsilon_{i+r}$, where $\theta_1 \theta_r \neq 0$. We shall also suppose that $r \geq 2$; the case $r = 1$ is straightforward. Let $J$ denote a large positive integer, and given $-\infty < j < \infty$, put $\mathcal{I}_j = (j/J, (j+1)/J]x$. Define $\varepsilon'_k = \theta_k \varepsilon_k$ and $\mathcal{E}_{jk} = \{\varepsilon'_k \in \mathcal{I}_j\}$, and let $\xi \in (0,1)$ be a constant. Suppose that the unique maximum of $\theta_k$ occurs at $k = \ell$. Then,

$P(X_0 > x$ and $\varepsilon_k > \xi x$ for some $k \in [1,r]$ with $k \neq \ell)$

$$\leq \sum_{\substack{i:\, 1 \leq i \leq r, \\ i \neq \ell}} P\left(\sum_{k=1}^r \varepsilon'_k > x, \varepsilon'_i > \xi x\right)$$

$$\leq \sum_{\substack{i:\, 1 \leq i \leq r, \\ i \neq \ell}} \sum_{\substack{j_1,\ldots,j_r: \\ j_1+\cdots+j_r+r \geq J}} P\left(\{\varepsilon'_i > \xi x\} \cap \bigcap_{k=1}^r \mathcal{E}_{j_k k}\right)$$

(5.8)
$$= \sum_{\substack{i:\, 1 \leq i \leq r, \\ i \neq \ell}} \sum_{\substack{j_1,\ldots,j_r: \\ j_1+\cdots+j_r+r \geq J}} P(\varepsilon'_i > \xi x, \mathcal{E}_{j_i i}) \prod_{\substack{k:\, 1 \leq k \leq r, \\ k \neq i}} P(\mathcal{E}_{j_k k})$$

$$\leq \exp\{o(x^\gamma)\} \sum_{\substack{i:\, 1 \leq i \leq r, \\ i \neq \ell}} \sum_{\substack{j_1,\ldots,j_r \geq 0 \\ j_1+\cdots+j_r+r \geq J}} \exp[-Cx^\gamma \max\{(\xi/\omega_2)^\gamma,$$

$$(j_i/J\theta_i)^\gamma\}]$$

$$\times \exp\left\{-Cx^\gamma \sum_{\substack{k:\, 1 \leq k \leq r, \\ k \neq i}} (j_k/J\theta_k)^\gamma\right\}.$$

If $i \neq \ell$, then the minimum of $\sum_{k:\, k \neq i}(u_k/\theta_k)^\gamma$, subject to $\sum_k u_k = v$ and each $u_k \geq 0$, occurs when $u_\ell = v$ and $u_k = 0$ for $k \neq \ell$. Hence, given $\xi > 0$, and $\eta > 0$ sufficiently small, we may choose $J$ so large that, uniformly in $1 \leq i \leq r$ with $i \neq \ell$,

$$\sum_{\substack{j_1,\ldots,j_r \geq 0: \\ j_1+\cdots+j_r+r \geq J}} \exp\Bigg[-Cx^\gamma \max\{(\xi/\omega_2)^\gamma, (j_i/J\theta_i)^\gamma\}$$

(5.9)
$$- Cx^\gamma \sum_{\substack{k:\, 1 \leq k \leq r, \\ k \neq i}} (j_k/J\theta_k)^\gamma\Bigg]$$

$$= O[\exp\{-(1+\eta)C(x/\omega_1)^\gamma\}].$$

Combining (5.8) and (5.9) we deduce that, for each $\xi \in (0,1)$, there exists $\eta = \eta(\xi) > 0$ for which, as $x \to \infty$,

$$P(X_0 > x \text{ and } \varepsilon_k > \xi x \text{ for some } k \in [1,r] \text{ with } k \neq \ell)$$

(5.10)
$$= O[\exp\{-(1+\eta)C(x/\omega_1)^\gamma\}] = O[\exp\{-C(x/\omega)^\gamma\}],$$

where $0 < \omega < \omega_1$.

Let $0 < \xi, \eta < 1$ and define $y = x - (r-2)\xi$ and $\omega_3 = \omega_1 + \omega_2$. Then, for each $i \neq 0$ and for $J$ sufficiently large, the argument leading to (5.9) gives

$$P(X_0 > x, X_i > x, \text{ and } \varepsilon_k \leq \xi x \text{ for all } k \in [1,r] \cup [1-i, r-i]$$
$$\text{except for } k = \ell \text{ or } k = \ell + i)$$

(5.11)
$$\leq P(\theta_\ell \varepsilon_\ell + \theta_{\ell+i} \varepsilon_{\ell+i} > y \text{ and } \theta_{\ell-i} \varepsilon_\ell + \theta_\ell \varepsilon_{\ell+i} > y)$$
$$\leq P\{(\theta_\ell + \theta_{\ell-i})\varepsilon_\ell + (\theta_\ell + \theta_{\ell+i})\varepsilon_{\ell+i} > 2y\}$$
$$= O\left(\sum_{\substack{j_1, j_2 \geq 0: \\ j_1 + j_2 + 2 > J}} \exp[-(1-\eta)C(2y/\omega_3)^\gamma\{(j_1/J)^\gamma + (j_2/J)^\gamma\}]\right)$$
$$= O[\exp\{-(1-\eta)C(2y/\omega_3)^\gamma\}] = O[\exp\{-C(x/\omega)^\gamma\}],$$

where $\omega$ can be taken in $(0, \omega_1)$ if $\eta$ is chosen sufficiently small. Combining (5.10) and (5.11) we deduce that, for some $0 < \omega < \omega_1$,

(5.12)
$$P(X_0 > x, X_j > x) = O[\exp\{-C(x/\omega)^\gamma\}].$$

Result (5.6), and hence part (i) of Theorem 3.4, follows from (5.7) and (5.12).

Next we derive part (ii) of Theorem 3.4. Let $\ell$ denote one of the $q$ distinct values of $k$ for which $\theta_k = \max\{\theta_1, \ldots, \theta_r\}$; write $\mathcal{I}$ for the set of indices $i$ such that $1 \leq |i| \leq r$; let $\mathcal{I}(\ell)$ be the set of $q-1$ indices $i \in \mathcal{I}$ which are such that $\theta_{\ell-i} = \theta_\ell$; and let $\mathcal{I}'(\ell)$ be the complement of $\mathcal{I}(\ell)$ in $\mathcal{I}$. Let $\bowtie_i$, for $i \in \mathcal{I}$, be a sequence composed of the inequalities $<$ or $>$, as in Section 3.7. Then, the probability $p(\mathcal{I})$ that $X_i \bowtie_i x$ for each $i \in \mathcal{I}$, and that, in addition, $X_0 > x$, is given by

$$p(\mathcal{I}) = \int I\Bigg\{\sum_{1 \leq k \leq r: \, k+i \neq \ell} \theta_k u_{k+i} + \theta_{\ell-i} u_\ell \bowtie_i x \text{ for } i \in \mathcal{I}(\ell),$$

$$\sum_{k=1}^{r} \theta_k u_{k+i} \bowtie_i x \text{ for } i \in \mathcal{I}'(\ell), \sum_{1 \leq k \leq r: \, k \neq \ell} \theta_k u_k + \theta_\ell u_\ell > x\Bigg\}$$

$$\times \left\{\prod_{k=-2r}^{2r} f(u_k)\right\} du_{-2r} \cdots du_{2r}.$$

Part (ii) of Theorem 3.4 can be derived by evaluating this integral, changing variable appropriately. The argument is outlined below.

Write $u'$ for the vector with components $u_{-2r}, \ldots, u_{2r}$, except that $u_\ell$ is excluded. Let $C$, $C_1$ and $C_2$ be as in (3.8). Then, changing variable from $u_\ell$ to $v = \theta_\ell u_\ell / x$, we have, as $x \to \infty$,

$$p(\mathcal{I}) = \frac{x}{\theta_\ell} \int I \left\{ \frac{1}{x} \sum_{1 \leq k \leq r\,:\, k+i \neq \ell} \theta_k u_{k+i} + v \bowtie_i 1 \text{ for each } i \in \mathcal{I}(\ell), \right.$$

$$\left. \frac{1}{x} \sum_{k=1}^{r} \theta_k u_{k+i} \bowtie_i 1 \text{ for each } i \in \mathcal{I}'(\ell), \frac{1}{x} \sum_{1 \leq k \leq r\,:\, k \neq \ell} \theta_k u_k + v > 1 \right\}$$

$$\times \left\{ \prod_{k \neq \ell} f(u_k) \right\} f\left( \frac{vx}{\theta_\ell} \right) du'\, dv$$

$$= \frac{x}{\theta_\ell} \int I\{ v \bowtie_i 1 \text{ for each } i \in \mathcal{I}(\ell), 0 \bowtie_i 1 \text{ for each } i \in \mathcal{I}'(\ell), v > 1 \}$$

$$\times \left\{ \prod_{k \neq \ell} f(u_k) \right\} f\left( \frac{vx}{\theta_\ell} \right) du'\, dv + o[x^{C_2+1-\gamma} \exp\{ -C(x/\theta_\ell)^\gamma \}]$$

$$= \frac{x}{\theta_\ell} I\{ \bowtie_i \text{ equals } > \text{ for each } i \in \mathcal{I}(\ell), \bowtie_i \text{ equals } < \text{ for each } i \in \mathcal{I}'(\ell) \}$$

$$\times \int_1^\infty f\left( \frac{vx}{\theta_\ell} \right) dv + o[x^{C_2+1-\gamma} \exp\{ -C(x/\theta_\ell)^\gamma \}]$$

$$= I\{ \bowtie_i \text{ equals } > \text{ for each } i \in \mathcal{I}(\ell), \bowtie_i \text{ equals } < \text{ for each } i \in \mathcal{I}'(\ell) \}$$

$$\times C^{-1} C_1 (x/\theta_\ell) x^{C_2+1-\gamma} \exp\{ -C(x/\theta_\ell)^\gamma \}$$

$$+ o[x^{C_2+1-\gamma} \exp\{ -C(x/\theta_\ell)^\gamma \}].$$

A similar but simpler argument shows that, as $x \to \infty$,

$$P(X_0 > x) \sim C^{-1} C_1 (x/\theta_\ell) x^{C_2+1-\gamma} \exp\{ -C(x/\theta_\ell)^\gamma \}.$$

Therefore, $p(\mathcal{I}) \sim P(X_0 > x)$ if "$\bowtie_i$ equals $>$ for each $i \in \mathcal{I}(\ell)$, and $\bowtie_i$ equals $<$ for each $i \in \mathcal{I}'(\ell)$"; while $p(\mathcal{I}) = o\{ P(X_0 > x) \}$ if the property in quotation marks fails. This result implies that $P(M = q, X_0 > x) \sim P(X_0 > x)$, which is equivalent to part (ii) of Theorem 3.4.

5.5. *Proof of Theorem 3.5.* Without loss of generality, $X_i = \theta_1 \varepsilon_{i+1} + \cdots + \theta_r \varepsilon_{i+r}$ for each $i$, where $\theta_1 \theta_r \neq 0$. Let $c_1 > 0$ be fixed but arbitrarily large, and define $\mathcal{I}_1 = (-\infty, -c_1]$, $\mathcal{I}_2 = (-c_1, c_1]$, $\mathcal{I}_3 = (c_1, x/r]$, and $\mathcal{I}_4 = (x/r, \infty)$. Put $\varepsilon_k' = \theta_k \varepsilon_k$ and $\mathcal{E}_{jk} = \{ \varepsilon_k' \in \mathcal{I}_j \}$, for $j = 1, \ldots, 4$. If none of $\varepsilon_1', \ldots, \varepsilon_r'$ is in $\mathcal{I}_4$, then $X_0 < x$. Moreover, the probability, $p(k, x)$ say, that just $k$ of $\varepsilon_1', \ldots, \varepsilon_r'$ are in $\mathcal{I}_4$, satisfies $p(k,x) \asymp x^{-k\rho}$ as $x \to \infty$. Therefore, if we define $\mathcal{E}_0 = \{X_0 > x\}$, $\mathcal{E}_4 = \{\text{exactly one of } \mathcal{E}_{41}, \ldots, \mathcal{E}_{4r} \text{ holds}\}$ and $\mathcal{E}_5 = \mathcal{E}_0 \cap \mathcal{E}_4$, then, as $x \to \infty$,

$$(5.13) \qquad\qquad P(\mathcal{E}_0 \setminus \mathcal{E}_5) = O(x^{-2\rho}).$$

Put $\mathcal{E}_{6i} = \mathcal{E}_{1i} \cup \mathcal{E}_{3i}$ and $\mathcal{E}_6 = \{\text{at least one of } \mathcal{E}_{61}, \ldots, \mathcal{E}_{6r} \text{ holds}\}$. Then,

$$P(\mathcal{E}_5 \cap \mathcal{E}_6) \leq \sum_{i_1 \neq i_2} \sum P(\mathcal{E}_{4i_1} \cap \mathcal{E}_{6i_2}) \leq \sum_{i_1 \neq i_2} \sum P(\theta_{i_1}\varepsilon > x/r) P(\theta_{i_2}|\varepsilon| > c_1)$$

$$\leq B_1 x^{-\rho} P(|\varepsilon| > c_1 \min_i \theta_i^{-1}),$$

where $B_1 > 0$ does not depend on $c_1$. Therefore,

$$\lim_{c_1 \to \infty} \limsup_{x \to \infty} x^\rho P(\mathcal{E}_5 \cap \mathcal{E}_6) = 0.$$

Combining this result and (5.13), and defining $\mathcal{E}_7 = \mathcal{E}_0 \cap \mathcal{E}_4 \cap \widetilde{\mathcal{E}}_6$, where $\widetilde{\mathcal{E}}_6$ denotes the complement of $\mathcal{E}_6$, we have,

$$(5.14) \qquad\qquad \lim_{c_1 \to \infty} \limsup_{x \to \infty} x^\rho P(\mathcal{E}_0 \setminus \mathcal{E}_7) = 0.$$

Let $c_2$ denote any fixed real number, and define $\mathcal{E}_{8i} = \{\varepsilon_j' \in (-c_1, c_1] \text{ for each } j \in [1, r] \text{ for which } j \neq i\}$, and

$$\mathcal{E}_9 = \mathcal{E}_9(c_1, c_2, x) = \bigcup_{i=1}^r \{\varepsilon_i' > x + c_2\} \cap \mathcal{E}_{8i}.$$

Since

$$\mathcal{E}_7 = \{X_0 > x\} \cap \{\text{exactly one of } \mathcal{E}_{41}, \ldots, \mathcal{E}_{4r} \text{ holds}\}$$
$$\cap \{\text{none of } \mathcal{E}_{61}, \ldots, \mathcal{E}_{6r} \text{ holds}\},$$

then

$$(5.15) \qquad \bigcup_{i=1}^r \{\varepsilon_i' > x + (r-1)c_1\} \cap \mathcal{E}_{8i} \subseteq \mathcal{E}_7 \subseteq \bigcup_{i=1}^r \{\varepsilon_i' > x - (r-1)c_1\} \cap \mathcal{E}_{8i}.$$

However, for each $i \in [1, r]$ and each $c_3 > 0$,

$$(5.16) \qquad\qquad P(x - c_3 \leq \varepsilon_i' \leq x + c_3) = o(x^{-\rho})$$

as $x \to \infty$. Writing $\triangle$ for the symmetric-difference binary operator, and combining (5.14)–(5.16), we deduce that

$$(5.17) \qquad\qquad \lim_{c_1 \to \infty} \limsup_{x \to \infty} x^\rho P(\mathcal{E}_0 \triangle \mathcal{E}_9) = 0.$$

Let $R_j(c_1, x)$, for $j \geq 1$, denote a function of $c_1$ and $x$ satisfying

$$\lim_{c_1 \to \infty} \limsup_{x \to \infty} x^\rho R_j(c_1, x) = 0,$$

let $\theta_{(1)} \geq \cdots \geq \theta_{(m)}$ denote a ranking of the $m$ nonzero $\theta_i$'s, define $\theta_{(q)} = 0$ for $q > m$ and put

$$\mathcal{E}(j) = \bigcup_{i=1}^{r} [\{\theta_i \varepsilon_{i+j} > x\}$$

$$\cap \{-c_1 < \theta_k \varepsilon_{k+j} \leq c_1 \text{ for each } k \in [1, r] \text{ for which } k \neq i\}].$$

Then $\mathcal{E}(0) = \mathcal{E}_9(c_1, 0, x)$, and so, by (5.17),

$$\begin{aligned}
P(M = q) \\
&= P(X_j > x \text{ for exactly } q \text{ values of } j \text{ satisfying } |j| \leq r) \\
&= P\{\mathcal{E}(j) \text{ holds for exactly } q \text{ values of } j \text{ satisfying } |j| \leq r\} \\
(5.18) \qquad &\quad + R_1(c_1, x) \\
&= \sum_{i=1}^{m} P\{\varepsilon > \theta_{i-j}^{-1} x \text{ for exactly } q \text{ values of } j \text{ satisfying } |j| \leq r\} \\
&\quad + R_2(c_1, x) \\
&= C x^{-\rho} m(\theta_{(q)}^{\rho} - \theta_{(q+1)}^{\rho}) + R_3(c_1, x).
\end{aligned}$$

Result (5.18) implies that $P(M_1 = q) \to p_q$, which is identical to $P(M_0 = q)$, completing the proof of Theorem 3.5.

5.6. *Proof of Theorem 3.6.* Let $U_1 = (X_{-r}, \ldots, X_{-1}, X_1, \ldots, X_r)^{\mathrm{T}}$ and $U_2 = X_0$, and define $U = (U_1^{\mathrm{T}}, U_2)^{\mathrm{T}}$, a $(2r+1)$-vector. Partition the covariance matrix, $\Sigma$, of $U$ in the ratio $2r : 1$, meaning that the top left-hand corner matrix, $\Sigma_{11}$ say, is $2r \times 2r$, the upper right-hand and lower left-hand matrices, $\Sigma_{12}$ and $\Sigma_{21}$, are $r \times 1$ and $1 \times r$, and the lower right-hand corner matrix is $1 \times 1$ and equals 1. In this notation, $U_1$, conditional on $U_2 = u$, is Normal $N(\Sigma_{12} u, \Sigma_{11} - \Sigma_{12}\Sigma_{21})$. In view of (3.16), $(\Sigma_{12})_{i1} = 1 - c_{|i|}\delta + o(\delta)$, $(\Sigma_{11})_{ij} = 1 - c_{|i-j|}\delta + o(\delta)$ and

$$\begin{aligned}
(\Sigma_{12}\Sigma_{21})_{ij} &= \{1 - c_{|i|}\delta + o(\delta)\}\{1 - c_{|j|}\delta + o(\delta)\} \\
&= 1 - \delta(c_{|i|} + c_{|j|}) + o(\delta),
\end{aligned}$$

and so $(\Sigma_{11} - \Sigma_{12}\Sigma_{21})_{ij} = \Sigma_1 \delta + o(\delta)$, where $(\Sigma_1)_{ij} = c_{|i|} + c_{|j|} - c_{|i-j|}$. Therefore, conditional on $U_2 = u$, $\delta^{-1/2}(U_1 - \Sigma_{12} u)$ is Normal $N(0, \Sigma_2)$, where $\Sigma_2 = \Sigma_1 + o(1)$ and does not depend on $u$. Hence, taking $Z = (Z_{-r}, \ldots, Z_{-1}, Z_1, \ldots, Z_r)$ to be Normal $N(0, \Sigma_2)$, we have:

$$P(X_i \bowtie_i t \text{ for } 1 \leq |i| \leq r | X_0 > t)$$

$$= \int_t^{\infty} P(X_i \bowtie_i t \text{ for } 1 \leq |i| \leq r | X_0 = u) \, du P(X_0 \leq u \mid X_0 > t)$$

$$= t \int_t^\infty P(X_i \bowtie_i t \text{ for } 1 \leq |i| \leq r | X_0 = u) \exp\left\{-\tfrac{1}{2}(u^2 - t^2)\right\} du + o(1)$$

$$= \int_t^\infty P(X_i \bowtie_i t \text{ for } 1 \leq |i| \leq r | X_0 = t + vt^{-1})$$

$$\times \exp(-v - \tfrac{1}{2}v^2 t^{-2}) \, dv + o(1)$$

$$= \int_0^\infty P[\{1 - c_{|i|}\delta + o(\delta)\}(t + vt^{-1})$$

$$+ \{1 + o(1)\}\delta^{1/2} Z_i \bowtie_i t \text{ for } 1 \leq |i| \leq r] e^{-v} \, dv + o(1)$$

$$= \int_0^\infty P(Z_i \bowtie_i \delta^{1/2} t c_{|i|} - \delta^{-1/2} t^{-1} v \text{ for } 1 \leq |i| \leq r) e^{-v} \, dv + o(1)$$

$$= \pi(\mathcal{S}) + o(1).$$

Adding over sequences $\mathcal{S}$ that include just $k$ ">" signs, we deduce that $\pi_k \to \pi_k^0$.

**Acknowledgments.** The authors are grateful for helpful comments of Serguei Novak and two reviewers.

## REFERENCES

BENJAMINI, Y. and HOCHBERG, Y. (1995). Controlling the false discovery rate: A practical and powerful approach to multiple testing. *J. Roy. Statist. Soc. Ser. B* **57** 289–300. MR1325392

BENJAMINI, Y. and HOCHBERG, Y. (2000). On the adaptive control of the false discovery fate in multiple testing with independent statistics. *J. Educ. Behav. Statist.* **25** 60–83.

BENJAMINI, Y. and YEKUTIELI, D. (2001). The control of the false discovery rate in multiple testing under dependency. *Ann. Statist.* **29** 1165–1188. MR1869245

BERNHARD, G., KLEIN, M. and HOMMEL, G. (2004). Global and multiple test procedures using ordered $p$-values—a review. *Statist. Papers* **45** 1–14. MR2019782

BLAIR, R. C., TROENDLE, J. F. and BECK, R.W. (1996). Control of familywise errors in multiple endpoint assessments via stepwise permutation tests. *Statist. Med.* **15** 1107–1121.

BROWN, B. W. and RUSSELL, K. (1997). Methods correcting for multiple testing: Operating characteristics. *Statist. Med.* **16** 2511–2528.

DUDOIT, S., SHAFFER, J. P. and BOLDRICK, J. C. (2003). Multiple hypothesis testing in microarray experiments. *Statist. Sci.* **18** 73–103. MR1997066

DUNNETT, C. W. and TAMHANE, A. C. (1995). Step-up testing of parameters with unequally correlated estimates. *Biometrics* **51** 217–227.

EFRON, B. (2007). Correlation and large-scale simultaneous significance testing. *J. Amer. Statist. Assoc.* **102** 93–103. MR2293302

FINNER, H. and ROTERS, M. (1998). Asymptotic comparison of step-down and step-up multiple test procedures based on exchangeable test statistics. *Ann. Statist.* **26** 505–524. MR1626043

FINNER, H. and ROTERS, M. (1999). Asymptotic comparison of the critical values of step-down and step-up multiple comparison procedures. *J. Statist. Plann. Inference* **79** 11–30. MR1704211

FINNER, H. and ROTERS, M. (2000). On the critical value behavior of multiple decision procedures. *Scand. J. Statist.* **27** 563–573. MR1795781

FINNER, H. and ROTERS, M. (2002). Multiple hypotheses testing and expected number of type I errors. *Ann. Statist.* **30** 220–238. MR1892662

GENOVESE, C. and WASSERMAN, L. (2004). A stochastic process approach to false discovery control. *Ann. Statist.* **32** 1035–1061. MR2065197

GODFREY, G. K. (1985). Comparing the means of several groups. *New Eng. J. Med.* **311** 1450–1456.

GOTZSCHE, P. C. (1989). Methodology and overt and hidden bias in reports of 196 double-blind trials of nonsteroidal antiinflammatory drugs in rheumatoid arthritis. *Control Clin. Trials* **10** 31–56.

HOCHBERG, Y. (1988). A sharper Bonferroni procedure for multiple tests of significance. *Biometrika* **75** 800–802. MR0995126

HOCHBERG, Y. and BENJAMINI, Y. (1990). More powerful procedures for multiple testing. *Statist. Med.* **9** 811–818.

HOCHBERG, Y. and TAMHANE, A. C. (1987). *Multiple Comparison Procedures.* Wiley, New York. MR0914493

HOLLAND, B. and CHEUNG, S. H. (2002). Familywise robustness criteria for multiple-comparison procedures. *J. Roy. Statist. Soc. Ser. B* **64** 63–77. MR1881845

HOMMEL, G. (1988). A comparison of two modified Bonferroni procedures. *Biometrika* **76** 624–625. MR1040659

KESSELMAN, H. J., CRIBBIE, R. and HOLLAND, B. (2002). Controlling the rate of Type I error over a large set of statistical tests. *Brit. J. Math. Statist. Psych.* **55** 27–39.

LEHMANN, E. L. and ROMANO, J. P. (2005). *Testing Statistical Hypotheses*, 3rd ed. Springer, New York. MR2135927

LEHMANN, E. L., ROMANO, J. P. and SHAFFER, J. P. (2005). On optimality of stepdown and stepup multiple test procedures. *Ann. Statist.* **33** 1084–1108. MR2195629

LUDBROOK, J. (1991). On making multiple comparisons in clinical and experimental pharmacology and physiology. *Clin. Exper. Pharm. Physiol.* **18** 379–392.

OLEJNIK, S., LI, J. M., SUPATTATHUM, S. and HUBERTY, C. J. (1997). Multiple testing and statistical power with modified Bonferroni procedures. *J. Educ. Behav. Statist.* **22** 389–406.

OTTENBACHER, K. J. (1991a). Statistical conclusion validity: An empirical analysis of multiplicity in mental retardation research. *Amer. J. Ment. Retard.* **95** 421–427.

OTTENBACHER, K. J. (1991b). Statistical conclusion validity—multiple inferences in rehabilitation research. *Amer. J. Phys. Med. Rehab.* **70** 317–322.

OTTENBACHER, K. J. (1998). Quantitative evaluation of multiplicity in epidemiology and public health research. *Amer. J. Epidemiology* **147** 615–619.

OTTENBACHER, K. J. and BARRETT, K. A. (1991). Measures of effect size in the reporting of rehabilitation research. *Amer. J. Phys. Med. Rehab.* **70** S131–S137.

PIGEOT, I. (2000). Basic concepts of multiple tests—A survey. *Statist. Papers* **41** 3–36. MR1746085

POCOCK, S. J., HUGHES, M. D. and LEE, R. J. (1987). Statistical problems in reporting of clinical trials. *J. Amer. Statist. Assoc.* **84** 381–392.

ROM, D. M. (1990). A sequentially rejective test procedure based on a modified Bonferroni inequality. *Biometrika* **77** 663–665. MR1087860

ROSENBERG, P. S., CHE, A. and CHEN, B. E. (2006). Multiple hypothesis testing strategies for genetic case-control association studies. *Statist. Med.* **25** 3134–3149. MR2247232

SARKAR, S. K. (1998). Some probability inequalities for ordered MTP2 random variables: A proof of the Simes conjecture. *Ann. Statist.* **26** 494–504. MR1626047

SARKAR, S. K. (2006). False discovery and false nondiscovery rates in single-step multiple testing procedures. *Ann. Statist.* **34** 394–415. MR2275247

SARKAR, S. K. and CHANG, C. K. (1997). The Simes method for multiple hypothesis testing with positively dependent test statistics. *J. Amer. Statist. Assoc.* **92** 1601–1608. MR1615269

SCHMIDT, R. and STADTMÜLLER, U. (2006). Nonparametric estimation of tail dependence. *Scand. J. Statist.* **33** 307–335. MR2279645

SCHMIDT, T. (2007). Coping with copulas. In *Copulas—From Theory to Applications in Finance* (J. Rank, ed.) 3–34. Risk Books, London.

SEN, P. K. (1999). Some remarks on Simes-type multiple tests of significance. *J. Statist. Plann. Inference* **82** 139–145. MR1736438

SHAO, Q.-M. (1999). A Cramér type large deviation result for Student's $t$-statistic. *J. Theoret. Probab.* **12** 385–398. MR1684750

SIMES, R. J. (1986). An improved Bonferroni procedure for multiple tests of significance. *Biometrika* **73** 751–754. MR0897872

SMITH, D. E., CLEMENS, J., CREDE, W., HARVEY, M. and GRACELY, E. J. (1987). Impact of multiple comparisons in randomized clinical trials. *Amer. J. Med.* **83** 545–550.

YEKUTIELI, D., REINER-BENAIM, A., BENJAMINI, Y., ELMER, G. I., KAFKAFI, N., LETWIN, N. E. and LEE, N. H. (2006). Approaches to multiplicity issues in complex research in microarray analysis. *Statist. Neerlandica* **60** 414–437. MR2291384

WANG, Q. (2005). Limit theorems for self-normalized large deviations. *Electron. J. Probab.* **10** 1260–1285. MR2176384

WRIGHT, S. P. (1992). Adjusted $p$-values for simultaneous inference. *Biometrics* **48** 1005–1013.

DEPARTMENT OF MATHEMATICS AND STATISTICS
UNIVERSITY OF MELBOURNE
MELBOURNE, VIC 3010
AUSTRALIA
E-MAIL: s.clarke@ms.unimelb.edu.au